\documentclass[10pt]{article}
\usepackage{amsmath}
\usepackage{amssymb}

\usepackage{cite}

\usepackage{hyperref}

\usepackage{lineno}

\usepackage{microtype}
\usepackage{graphicx}
\DisableLigatures[f]{encoding = *, family = * }

\usepackage{algorithmic}
\usepackage{algorithm}
\usepackage[labelfont=bf,labelsep=period,justification=raggedright]{caption}

\makeatletter
\renewcommand{\@biblabel}[1]{\quad#1.}

\newcommand{\thma}{\begin{thm}}
\newcommand{\thmb}{\end{thm}}

\newcommand{\mata}{\begin{bmatrix}}
\newcommand{\matb}{\end{bmatrix}}

\newtheorem{thm}{Theorem}

\newtheorem{defi}{Definition}

\newtheorem{prop}{Proposition}

\newcommand{\enuma}{\begin{enumerate}}
\newcommand{\enumb}{\end{enumerate}}
\newcommand{\ena}{\begin{enumerate}}
\newcommand{\enb}{\end{enumerate}}

\newcommand{\itema}{\begin{itemize}}
\newcommand{\itemb}{\end{itemize}}
\newcommand{\ita}{\begin{itemize}}
\newcommand{\itb}{\end{itemize}}

\newcommand{\proofa}{\begin{proof}}
\newcommand{\proofb}{\end{proof}}

\newcommand{\bla}{\begin{block}}
\newcommand{\blb}{\end{block}}

\providecommand{\mc}[1]{\mathcal{#1}}
\providecommand{\mb}[1]{\boldsymbol{#1}}

\providecommand{\mh}[1]{\hat{#1}}
\providecommand{\wh}[1]{\widehat{#1}}

\providecommand{\mt}[1]{\widetilde{#1}}

\newcommand{\PmcP}{P \in \mc{P}}

\newcommand{\defa}{\begin{defi}}
\newcommand{\defb}{\end{defi}}

\providecommand{\mc}[1]{\mathcal{#1}}
\providecommand{\mb}[1]{\boldsymbol{#1}}

\providecommand{\mh}[1]{\hat{#1}}
\providecommand{\wh}[1]{\widehat{#1}}

\providecommand{\mt}[1]{\widetilde{#1}}

\newcommand{\T}{^{\ensuremath{\mathsf{T}}}}

\newcommand{\argmin}{\operatornamewithlimits{argmin}}

\providecommand{\ve}[1]{\boldsymbol{#1}}

\newtheorem{rem}{Remark}

\makeatother

\date{}

\begin{document}
% Title must be 150 characters or less
\begin{flushleft}
{\Large
\textbf{Fast Approximate Quadratic Programming for\\ Large (Brain) Graph Matching}
}
% Insert Author names, affiliations and corresponding author email.
\\
Joshua T. Vogelstein$^{1,\ast}$, John M. Conroy$^{2}$, Vince Lyzinski$^{3}$,
, Louis J. Podrazik$^{2}$,  Steven G.~Kratzer$^{2}$, Eric T.~Harley$^{4}$,
Donniell E.~Fishkind$^{4}$, 
		R.~Jacob~Vogelstein$^{5}$
        and Carey E.~Priebe$^{4}$
\\        
\bf{1} Department of Biomedical Engineering, Johns Hopkins University, Baltimore, MD, USA
\\
\bf{2} Institute for Defense Analyses, Center for Computing Sciences, Bowie, MD, USA
\\
\bf{3} Johns Hopkins University Human Language Technology Center of Excellence, Baltimore, MD, USA
\\
\bf{4} Department of Applied Mathematics and Statistics, Johns Hopkins University, Baltimore, MD, USA
\\
\bf{5}  Johns Hopkins University Applied Physics Laboratory, Laurel, MD, USA
\\
$\ast$ E-mail: jovo@jhu.edu
\end{flushleft}

\section*{Abstract}
Quadratic assignment problems arise in a wide variety of domains, spanning operations research, graph theory, computer vision, and neuroscience, to name a few.  The graph matching problem is a special case of the quadratic assignment problem,
and graph matching is increasingly important as
graph-valued data is becoming more prominent.  With the aim of efficiently and accurately matching the large graphs common in big data,
we present our graph matching algorithm, the Fast Approximate Quadratic assignment algorithm. 
We empirically demonstrate that our algorithm is faster and achieves a lower objective value on over $80\%$ of the QAPLIB benchmark library, compared with the previous state-of-the-art.  Applying our algorithm to our motivating example, matching C. elegans connectomes (brain-graphs), we find that it efficiently achieves  performance.

\section{Introduction}

In its most general form, the graph matching problem (GMP)---finding an alignment of the vertices of two graphs which minimizes the number of induced edge disagreements---is equivalent to a quadratic assignment problem (QAP) \cite{Umeyama1988}.
QAPs were first devised by Koopmans and Beckmann in 1957 to solve a ubiquitous problem in distributed resource allocation \cite{Koopmans1957}, and many important problems in combinatorial optimization (for example, the ``traveling salesman problem,'' and the GMP) have been shown to be specialized QAPs.  
While QAPs are known to be {\bf NP}-hard in general \cite{Papadimitriou1998},
they are widely applicable and there is a large literature devoted to their approximation and formulation; see \cite{burkard} for a comprehensive literature survey.
In casting the GMP as a QAP, we bring to bear a host of existing optimization theoretic tools and methodologies for addressing graph matching: Frank-Wolfe \cite{Frank1956}, gradient-descent \cite{bazaraa2013nonlinear}, etc.

Graph matching has applications in a wide variety of disciplines, spanning computer vision, image analysis, pattern recognition, and neuroscience; see \cite{Conte2004} for a comprehensive survey of the graph matching literature.  
We are motivated by applications in ``connectomics,''  an emerging discipline within neuroscience devoted to the study of brain-graphs, where vertices represent (collections of) neurons and edges represent connections between them \cite{Sporns05a, Hagmann05a}.  Analyzing connectomes is an important step for many neurobiological inference tasks.  For example, it is becoming increasingly popular to diagnose neurological diseases via comparing brain images \cite{Csernansky2004}.  To date, however, these comparisons have largely rested on anatomical (e.g., shape) comparisons, not graph (e.g., structural) comparisons.  This is despite the widely held doctrine that many
% Yet almost immediately after the ``neuron doctrine'' was conjectured (the idea that networks of neurons comprise brains), Wernicke and others began postulating that 
psychiatric disorders are fundamentally ``connectopathies,'' i.e. disorders of the connections of the brain \cite{Kubicki2007,Calhoun2011,Fornito2012,Fornito2012a}.   
Thus, available tests for connectopic explanations of psychiatric disorders hinge on first choosing particular graph invariants to compare across populations, rather than comparing the graphs' structure directly. Yet, recent results suggest that the graph invariant approach to classifying is both theoretically and practically inferior to comparing whole graphs via matching \cite{VP11_unlabeled}.

Part of the reason for the lack of publications that structurally compare brain-graphs is that existing algorithms for matching very large graphs are largely ineffective, often sacrificing accuracy for efficiency.
Indeed, available human connectomes have $\mc{O}(10^6)$ vertices and $\mc{O}(10^8)$ edges, and building them leverages state-of-the-art advances in DT-MRI imagery, big data processing and computer vision \cite{MIGRAINE}.  %Typically, in human connectomics, the brain is subdivided into approximately 100 vertices (regions), 
Contrast this with the fact that the human brain consists of approximately 86 billion neurons \cite{Herculano-Houzel2012}.   
At the other end of the spectrum, the small hermaphroditic  \emph{Caenorhabditis elegans} worm (\emph{C. elegans}) has only 302 neurons, with a fully mapped connectome.  
%Thus, regardless of which brains are under investigation, they are currently typically represented by graphs with $\mc{O}(100)$ vertices.
Consequently, these graphs demand GM algorithms that both accurately match small--to--moderately sized graphs and also
%these applications demand algorithms that 
scale to match very large graphs. 

% Moreover, 
% state-of-the-art procedures for many decision-theoretic or statistical inference tasks follow from first constructing interpoint dissimilarity matrices, i.e. graphs, from the data \cite{Duin2011}.  In this paradigm, graph matching is becoming a fundamental subroutine of many statistical inference pipelines. 

% As the number of vertices increases, we are forced to use approximate algorithms to find good (but not optimal) solutions.  

When matching these large connectome graphs (and, more broadly, the large graphs common in big data \cite{Kolaczyk2010}), we
necessarily 
% exact matching is often intractable.   Instead, we employ inexact matching algorithms (also called ``heuristics'') that will scale polynomially or even linear \cite{Conte2004}.
% For these problems, 
face an accuracy/efficiency trade-off when approximately solving these GMPs: slower algorithms could achieve better performance given more time (at the extreme, exhaustive search algorithms clearly have optimal performance).
%  Thus for any applied problem, 
%given the no free lunch theorem, 
Our algorithm---the Fast Approximate QAP (FAQ) algorithm for approximate GM---achieves the best available trade-off between accuracy and efficiency, outperforming the existing state-of-the-art in both accuracy and efficiency on a large proportion of QAP benchmarks and biologically inspired network matching problems.

The remainder of this paper is organized as follows. 
Section \ref{sec:QAP} formally defines the QAP and a relaxation thereof that we will operate under.  Section \ref{sec:GM} defines graph matching, and explains how it can be posed as a QAP. 

Section \ref{sec:FAQ} describes the FAQ algorithm.  Section \ref{sec:results} provides a number of theoretical and empirical results, and compares our algorithm to previous state-of-the-art algorithms.  This section concludes with an analysis of FAQ on our motivating problem of matching worm brain connectomes. We conclude with a discussion of possible extensions of FAQ and related work in Section \ref{sec:discussion}.

\section{Preliminaries}
In this section, we formally define the QAP and the GMP.  We then show how the GMP can be recast as a special case of the QAP.
\subsection{Quadratic Assignment Problems} 
\label{sec:QAP}

We first define the general quadratic assignment problem.  
Let $A=(a_{uv}),\,B=(b_{uv})\in\mathbb{R}^{n\times n}$ be two $n \times n$ real matrices.
Let $\Pi:=\Pi_n$ be the set of permutation functions (bijections) of the set
$[n]=\{1,\ldots, n\}$. We define the Koopmans-Beckmann (KB) version QAP via:

\begin{equation}
\text{(KB)} \qquad  
\begin{array}{cl}
			\text{minimize}   &\sum_{u,v \in [n]} b_{uv}a_{\pi(u)\pi(v)} \\
			\text{subject to}  &\pi \in \Pi  .
\end{array}\label{eq:KB}
\end{equation}
Note that occasionally an additional linear function is added, though we drop it here for brevity.

Equation (\ref{eq:KB}) can also be recast in matrix notation. Let $\mc{P}$ be the set of  $n \times n$ \emph{permutation matrices}, $\mc{P}=\{P\in \{0,1\}^{n \times n} : P\T \mb{1} = P \mb{1} = \mb{1} \}$, where $\mb{1}$ is the $n$-dimensional column vector consisting of all $1$'s. 
% Thus, we can write $PAP\T=(a_{\pi(u)\pi(v)})$ whenever $P$ is the permutation matrix corresponding to the bijection $\pi$, yielding the following equivalent optimization problem:
% \begin{equation*}
% % \text{(QAP)} \qquad  
% \begin{array}{cl}
% 			\text{minimize}   &\sum_{u,v \in [n]} b_{uv} (p_{vu} a_{uv} p_{uv}) \\
% 			\text{subject to}  &P \in \mc{P}.   
% \end{array}\label{eq:QAP}
% \end{equation*}
Thus, Equation (\ref{eq:KB}) can be written more compactly in matrix notation via:

\begin{equation}
\text{(QAP)} \qquad  
\begin{array}{cl}
			\text{minimize}   & \text{trace}(APB\T P\T) \\
			\text{subject to}  &P \in \mc{P}.   
\end{array}\label{eq:QAP}
\end{equation}
We hereafter refer to (\ref{eq:QAP}) as the QAP optimization function.

\subsection{Relaxed Quadratic Assignment Problem}

Eq. (\ref{eq:QAP}) is a binary quadratic program with linear constraints. 
%(note that it can also be written as a quadratic problem with quadratic constraints, because $p_{uv} \in \{0,1\}$ is equivalent to $p_{uv}=p_{uv}^2$).  Thus, one could use any number of algorithms for solving quadratic problems with binary (or quadratic) constraints.  
Because of the combinatorial nature of the feasible region, finding a \emph{global} optimum of (\ref{eq:QAP}) is \textbf{NP}-hard.  
Rather than directly optimizing over the permutation matrices, we begin by relaxing the constraint set to the convex hull of $\mc{P}$, the set of doubly stochastic matrices (i.e. the Birkhoff polytope),

$$\mc{D}:=\mc{D}_n=\{P\in \mathbb{R}^{n\times n} : P\T \mb{1} =  P \mb{1} = \mb{1}, P\succeq 0 \},$$ where $\succeq$ indicates an element-wise inequality. Relaxing $\mathcal{P}$ to $\mc{D}$ in (\ref{eq:QAP}) yields the relaxed quadratic assignment problem (rQAP):

\begin{equation}
\text{(rQAP)} \qquad  
\begin{array}{cl}
			\text{minimize}   & \text{trace}(APB\T P\T) \\
			\text{subject to}  &P \in \mc{D}.   
\end{array}\label{eq:rQAP}
\end{equation}
Note that, although rQAP is a quadratic program with linear constraints, it is {\it not} necessarily convex.  Indeed, the objective function, $f(P)= \text{trace}(APB\T P\T)$, has a Hessian that is not necessarily positive definite:
	$\nabla^2 f(P)  =   B \otimes A + B\T \otimes A\T,$
where $\otimes$ indicates the Kronecker product (note that if $A$ and $B$ are hollow matrices---as is common for graphs---then $\nabla^2 f(P)$ has trace equal to 0, and is indefinite).

While nonconvex quadratic optimization is, in general, {\bf NP}-hard, relaxing the feasible region allows us to employ the tools of continuous optimization to
search for a \emph{local} optima of (\ref{eq:rQAP}).  These local optima can then be projected onto $\mc{P}$, yielding an approximate solution of (\ref{eq:QAP}).
We also note that when relaxed to $\mc{D}$, the QAP optimization function is often multimodal, making initialization important when solving (\ref{eq:rQAP}).  
%With this in mind, below, we describe an algorithm to find a local optimum of rQAP.

\subsection{Graph Matching} % (fold)
\label{sec:GM}

A labeled graph $G=(\mc{V},\mc{E})$ consists of a vertex set $\mc{V}=[n]$, and an edge set $\mc{E}\subset\binom{\mc{V}}{2}$ in the undirected case, or 
$\mc{E}\subset\mc{V}\times\mc{V}$ in the directed case. %, where $|\mc{E}| \leq n^2$.
For an $n$-vertex graph $G$, we define the associated adjacency matrix $A=(a_{uv})\in\{0,1\}^{n\times n}$ to be the binary $n\times n$ matrix with
$a_{uv}= 1 \text{ if }\{u,v\}\in \mc{E}\text{ in the undirected setting},$ or $(u,v)\in \mc{E}\text{ in the directed setting}.$
%Note that we are not restricting our formulation to be directed or exclude self-loops. 
Given a pair of $n$-vertex graphs $G_A=(\mc{V}_A,\mc{E}_A)$ and $G_B=(\mc{V}_B,\mc{E}_B)$, with respective adjacency matrices $A$ and $B$, 
% let $\Pi$ be the set of permutation functions (bijections), $\Pi=\{\pi \from \mc{V}_A \to \mc{V}_B\}$.
% $\pi: \mc{V}_1 \to \mc{V}_2$ be a permutation function (bijection), and let $\Pi$ be the set of all such permutation functions.  
% Now 
we consider the following two closely related problems:
% A pair of graphs, $G_1$ and $G_2$, are isomorphic if and only if the following \emph{isomorphism criterion} holds: there exists a $\pi \in \Pi$ such that . 
% Let $A$ be the adjacency matrix representation of graph such that $A_{ij}=1$ if there is an edge from $u$ to $v$, and $A_{ij}=0$ otherwise. 
% Note that the below follows for directed/undirected and loopy/non-loopy graphs.
% $u \sim v \in \mc{E}$ and $A_{ij}=0$ otherwise.  
% Let  $\Pi$ be the set of permutation functions, where a permutation function (bijection) $\pi: \mc{V} \to \mc{V}$ (re-)orders the elements of the set $\mc{V}$.  Given a pair of $n \times n$ adjacency matrices, $A=(a_{ij})$ and $B=(b_{ij})$, consider the following two problems:
\begin{itemize}
	\item \textbf{Graph Isomorphism (GI):}  Does there exist a $P \in \mc{P}$ such that $A=PBP\T$. 
		\item \textbf{Graph Matching:}
		% Which $\pi \in \Pi$ minimizes the number of pairs of vertices $u,v \in \mc{V}_A$ such that $(u,v) \in \mc{E}_A$ and $(\pi (u) ,\pi (v)) \not \in \mc{E}_B$ or $(u,v) \not \in \mc{E}_A$ and  $(\pi (u) ,\pi (v)) \in \mc{E}_B$
		 $\min_{P\in\mc{P}}\|A-PBP\T \|_F$, where $\|\cdot\|_F$ is the usual matrix Frobenius norm.
\end{itemize}
GI is one of few problems with unknown computational complexity in the {\bf NP}-hierarchy \cite{Fortin1996};  indeed, if \textbf{P}$\neq$\textbf{NP}, then GI might reside in an intermediate complexity class called \textbf{GI}-complete.  Moreover, GI is, at worst, only moderately exponential, with complexity $\mc{O}(\exp\{n^{1/2 + o(1)}\})$ \cite{Babai1981}. On the other hand, the (harder) GMP---reducible to a QAP---is known to be \textbf{NP}-hard in general.    
Although polynomial time algorithms are available for GM (and GI) for large classes of problems (e.g., planar graphs, trees) \cite{Babai1980}, these
% There exist no known algorithms for which worst case behavior is polynomial \cite{Fortin1996}.  While GM is known to be \textbf{NP}-hard, it remains unclear whether GI is in $\mc{P}$, \textbf{NP}, or its own intermediate complexity class, \textbf{NP}-isomorphism (or isomorphism-complete).  
algorithms often have lead constants which are very large.  For example, if all vertices have degree less than $k$, there is a linear time algorithm for GI.  However, the hidden constant in this algorithm is $(512k^3)!$ \cite{Chen1994}.
Because we are interested in solving GM for graphs with $\dot{\approx} 10^6$ or more vertices, exact GM solutions will be computationally intractable. As such, we develop a fast approximate graph matching algorithm.   Our approach is based on formulating GM as a quadratic assignment problem.

\subsection{Graph Matching as a Quadratic Assignment Problem}
Given a pair of $n$-vertex graphs $G_A=(\mc{V}_A,\mc{E}_A)$ and $G_B=(\mc{V}_B,\mc{E}_B)$, with respective adjacency matrices $A$ and $B$, 
we can formally write the graph matching problem as an optimization problem:

\begin{equation}
% \text{(GM)} \qquad  
\begin{array}{cl}
			\text{minimize}   &\|AP-PB\|^2_F \\
			\text{subject to}  &P\in\mc{P}.   
\end{array}\label{eq:GMc}
\end{equation}
Simple algebra yields that,
\begin{align} \label{eq:equiv}
\|AP - PB\|_F^2  &= \text{trace} \{ (AP - PB)\T (AP - PB)\}\notag\\
&=\text{trace}(A\T A)+\text{trace}(BB\T)-2\text{trace}(APB\T P\T).
\end{align}
The GMP is then equivalent (i.e. same $\argmin$) to 

\begin{equation}
\text{(GM)} \qquad  
\begin{array}{cl}
			\text{minimize}   & - \text{trace}(APB\T P\T) \\
			\text{subject to}  &P \in \mc{P}.   
\end{array}\label{eq:GM}
\end{equation}
The objective function for GM is just the negative of the objective function for QAP. Thus, any descent algorithm for the former can be directly applied to the latter.  Moreover, any QAP approximation algorithms also immediately yields an analogous GM approximation.

As is common in solving general QAPs, GM algorithms often begin by first relaxing (\ref{eq:GMc}) to a continuous problem (see, for example, \cite{Zaslavskiy2009}).  The resulting problem is a convex quadratic program and can be efficiently exactly solved.  The obtained solution is then projected back onto $\mc{P}$ yielding an approximate solution to (\ref{eq:GMc}).  Contrary to popular existing approaches, our FAQ algorithm first solves a relaxed version of (\ref{eq:GM}) and subsequently projects the solution back onto $\mc{P}$.  
This relaxation yields an indefinite quadratic program, and indefinite quadratic programs are are {\bf NP}-hard to solve in general.  
However, recent theory indicates that the indefinite relaxation of (\ref{eq:GM}), and {\it not} the convex relaxation of (\ref{eq:GMc}) is the provably correct approach \cite{lyzinski2014graph}.
Reflecting this theory, we find that FAQ obtains state-of-the-art performance in terms of both computational efficiency and objective function value for various QAPs (see Section \ref{sec:results}).

\section{Fast Approximate QAP Algorithm} 
\label{sec:FAQ}

Our algorithm, called FAQ, proceeds in three steps:
\begin{enumerate}
	\item[] A. Choose a suitable initial position
	\item[] B. Find a local solution to rQAP. 
	\item[] C. Project back onto the set of permutation matrices. 
\end{enumerate}

These steps are summarized in Algorithm \ref{alg:1}.  Below, we provide details for each step.

\begin{algorithm}
	\caption*{FAQ for finding a local optimum of rQAP} 
\begin{algorithmic}[1]
\caption{FAQ for finding a local optimum of rQAP} \label{alg:1}
	\REQUIRE Graphs (adjacency matrices) $A$ and $B$ as well as a stopping criteria
	\ENSURE $\wh{P}$, an estimated permutation matrix
	\STATE Choose an initialization, $P^{(0)}=\mb{1}\mb{1}\T/n$ \label{step:init} 
	\WHILE{stopping criteria not met} 
	\STATE Compute the gradient of $f$ at the current point via Eq. \eqref{eq:grad}	
	\STATE Compute the direction $Q^{(i)}$ by solving Eq. \eqref{eq:dir} via the Hungarian algorithm
	\STATE Compute the step size $\alpha^{(i)}$ by solving Eq. \eqref{eq:step}
	\STATE Update $P^{(i)}$ according to Eq. \eqref{eq:update}  %$P^{(i+1)} \leftarrow P^{(i)} + \alpha^{(i)} Q^{(i)}$
	\ENDWHILE
	\STATE Obtain $\wh{P}$ by solving Eq. \eqref{eq:proj} via the Hungarian algorithm.
\end{algorithmic}
\end{algorithm}

\noindent\textbf{A: Find a suitable initial position.}  While any doubly stochastic matrix would be a feasible initial point, we choose the 
noninformative
``flat doubly  stochastic matrix,'' $J=\ve{1} \cdot \ve{1}\T/n$, i.e. the barycenter of the feasible region.
Alternately, we could use multiple restarts, each initial point near $J$.  Specifically, we could sample $K$, a random doubly stochastic matrix using 10 iterations of Sinkhorn balancing \cite{Sinkhorn1964}, and let $P^{(0)}=(J+K)/2$. Given this initial estimate (or estimates), we would then iterate the following five steps until convergence.

\noindent\textbf{B: Find a local solution to rQAP.} As mentioned above, rQAP is a quadratic problem with linear constraints.  A number of off-the-shelf algorithms are readily available for finding local optima in such problems.  We utilize the Frank-Wolfe algorithm (FW), a successive first-order optimization procedure originally devised to solve convex quadratic programs \cite{Frank1956, Bradley1977}.
Although FW is a relatively standard solver, especially as a subroutine for QAP algorithms \cite{Anstreicher03}, we provide a detailed view of applying FW to rQAP.

Given an initial position, $P^{(0)}$, iterate the following four steps:\\
{\it Step 1,} {\it Compute the gradient $\nabla f(P^{(i)})$:}  The gradient of $f(P)=-\text{trace}(APB\T P\T)$ with respect to $P$, evaluated at $P^{(i)}$, is given by $\label{eq:grad}
	\nabla f (P^{(i)}) = 
	% \partial f / \partial P^{(i)} =
	  - A P^{(i)} B\T - A\T P^{(i)} B.$\\
{\it Step 2,} {\it Compute the search direction $Q^{(i)}$:} The search direction is given by the argument that minimizes a first-order Taylor series approximation to $f(P)$ around the current estimate, $P^{(i)}$:

\begin{equation}
	\mt{f}^{(i)}(P) := f(P^{(i)}) + \nabla f(P^{(i)})\T(P - P^{(i)}).
\end{equation}
Dropping terms independent of $P$, we obtain the following sub-problem:

\begin{equation}
% \text{(rQAP)} \qquad  
\begin{array}{cl}
			\text{minimize}   & \text{trace}(\nabla f(P^{(i)})\T P) \\
			\text{subject to}  &P \in \mc{D},   
\end{array} \label{eq:dir}
\end{equation}
which is equivalent to a \emph{Linear Assignment Problem} (LAP), and can therefore be solved in $O(n^3)$ time via the ``Hungarian Algorithm'' of \cite{Kuhn1955, Jonker1987}.  
Let $Q^{(i)}$ indicate the argmin of Eq. \eqref{eq:dir}.
\\
{\it Step 3,} {\it Compute the step size $\alpha^{(i)}$:} Given $Q^{(i)}$, we minimize the original optimization problem, along the line segment from $P^{(i)}$ to $Q^{(i)}$:    

\begin{equation}
% \text{(rQAP)} \qquad  
\begin{array}{cl}
			\text{minimize}   & f(P^{(i)} + \alpha^{(i)} Q^{(i)}) \\
			\text{subject to}  & \alpha \in [0,1].   
\end{array} \label{eq:step}
\end{equation}
This can be solved exactly, as $f$ is a quadratic function of $\alpha$.  Let $\alpha^{(i)}$ indicate the $\argmin$ of Eq. \eqref{eq:step}.
\\
{\it Step 4,} {\it Update $P^{(i)}$:} Finally, the next iterate is the doubly stochastic matrix 

\begin{equation} \label{eq:update}
	P^{(i+1)} = P^{(i)} + \alpha^{(i)} Q^{(i)}.
\end{equation}
\\
{\it Stopping criteria:} Steps 1--4 are iterated until some stopping criterion is met.  Often, a threshold $\epsilon>0$ or an iteration limit is given, and the algorithm iterates until the iteration limit is reached, $\|P^{(i)}-P^{(i-1)}\|_F<\epsilon$, or $\|\nabla f(P^{(i)})\|_F<\epsilon$.  In practice, the algorithm often converges with a modest number of FW iterates.

\noindent\textbf{C: Project onto the set of permutation matrices.}   Let $P^{(final)}$ be the doubly stochastic matrix resulting from the final iteration of FW.  We project $P^{(final)}$ onto the set of permutation matrices via

\begin{equation}
% \text{(rQAP)} \qquad  
\begin{array}{cl}
			\text{minimize}   & -\text{trace}( P^{(final)} P\T) \\
			\text{subject to}  & \PmcP.   
\end{array} \label{eq:proj}
\end{equation}
  Note that Eq. \eqref{eq:proj} is again equivalent to a LAP and can be solved in $O(n^3)$ time.

\section{Results} 
\label{sec:results}

Below we present a number of empirical and theoretical results demonstrating the state-of-the-art efficiency and accuracy of the FAQ algorithm.

\subsection{Algorithm Complexity and leading constants} % (fold)
\label{sub:const}

\begin{figure}[htbp]
	\centering			
	\includegraphics[width=0.7\linewidth]{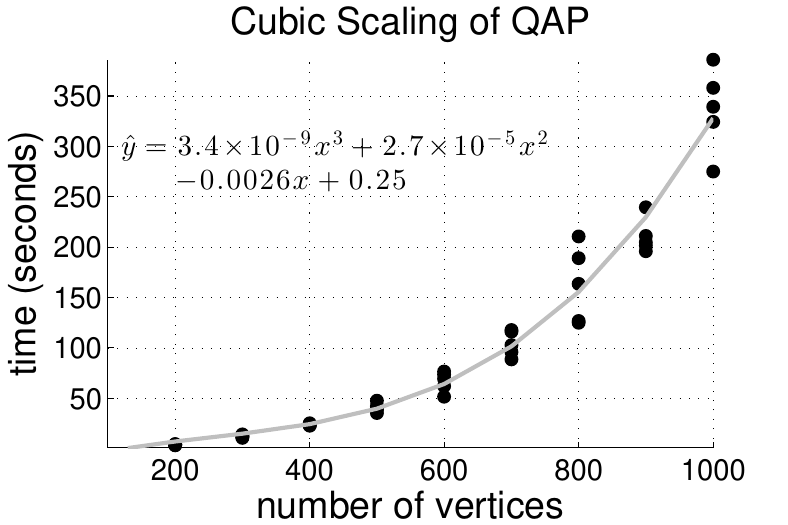}
	\caption{{\bf Running time of FAQ as function of number of vertices}. Data was sampled from an Erd\"os-R\'enyi model with $p=log(n)/n$.  Each dot represents a single simulation, with 100 simulations per $n$.  The solid line is the best fit cubic function.  Note the leading constant is $\approx 10^{-9}$. FAQ finds the optimal objective function value in every simulation.}
	\label{fig:scaling}
\end{figure}

As mentioned above, GM is computationally difficult, and even in the special cases for which polynomial time algorithms are available, the leading constants are intractably large. Given a bounded number of FW iterates, the FAQ algorithm has complexity $O(n^3)$.  However, a very large lead order constant could still render this algorithm practically infeasible.  
Figure \ref{fig:scaling} suggests that FAQ has $O(n^3)$ complexity, and also has very small leading constants ($\approx 10^{-9}$). 
This suggests that this algorithm is feasible for matching even reasonably large graphs.  Note that other state-of-the-art approximate graph matching algorithms also have cubic or worse time complexity in the number of vertices.  We will describe these other algorithms and their time complexity in greater detail below.

\subsection{QAP Benchmark Accuracy} 
\label{sub:qap_benchmarks}

Having demonstrated both theoretically and empirically that FAQ has cubic time complexity, we evaluate the algorithms accuracy on a suite of standard benchmarks.  More specifically, QAPLIB is a library of 137 quadratic assignment problems, ranging in size from 10 to 256 vertices \cite{Burkard1997}.  Recent graph matching papers typically evaluate the performance of their algorithm on 16 of the benchmark QAPs that are known to be particularly difficult \cite{Zaslavskiy2009,Schellewald2001}.  We compare the results of FAQ to the results of four other state-of-the-art graph matching algorithms: 
\begin{itemize}
\item[](1) the PATH algorithm, which solves the relaxation of (\ref{eq:GMc}), and then solves a sequence of concave and convex problems to project the solution onto $\mc{P}$ \cite{Zaslavskiy2009};
\item[](2) QCV which solves the relaxation of (\ref{eq:GMc}), and projects the obtained solution onto the closest permutation in $\mc{P}$;
\item[](3) the RANK algorithm \cite{Singh2007}, a spectral graph matching procedure;
\item[](4) the Umeyama algorithm (denoted by U), another spectral graph matching procedure \cite{Umeyama1988}.
\end{itemize}
  Note that the code for these four algorithms is freely available from the graphm package \cite{Zaslavskiy2009}.  

\begin{figure}[htbp]
	\centering
		\includegraphics[width=1.0\linewidth]{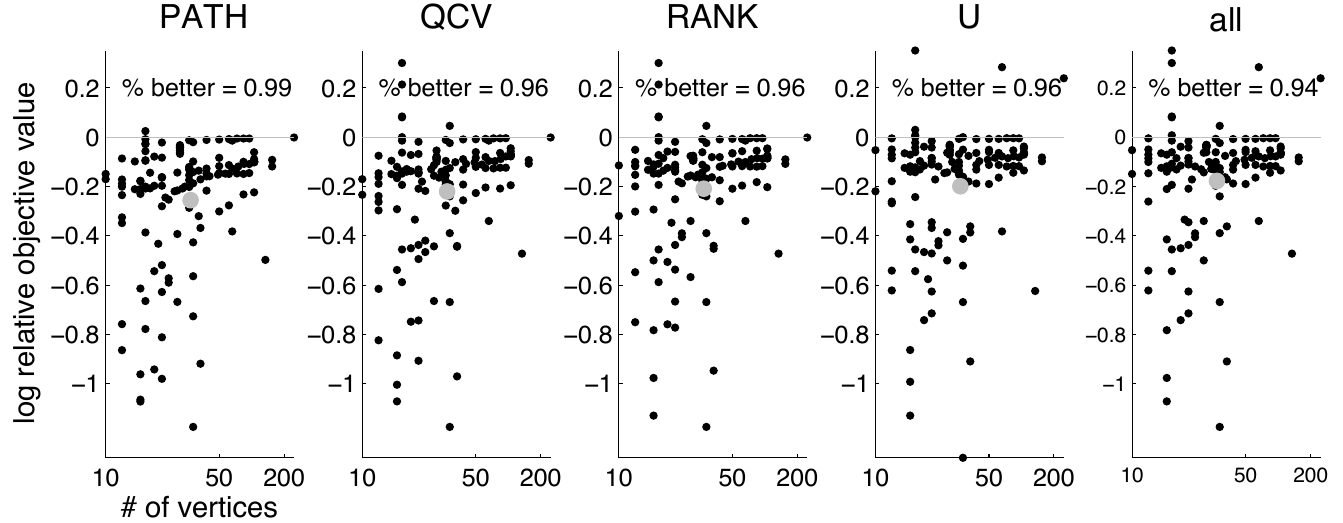}
	\caption{{\bf Relative accuracy---defined to be $\log_{10}(\mh{f}_{FAQ}/\mh{f}_X)$---of all the four algorithms compared with FAQ}.  Note that FAQ is better than all the other algorithms on $\approx 94\%$ of the benchmarks. The abscissa is the log number of vertices.  The gray dot indicates the mean improvement of FAQ over the other algorithms.}
	\label{fig:allRelAccuracy}
\end{figure}

Performance is measured by minimizing the assignment cost $f(P)=\text{trace}(APB\T P\T)$.  We write $\mh{f}_X$ for the value of the local minimum of $f$ obtained by the algorithm $X \in \{$FAQ, PATH, QCV, RANK, U, all$\}$, where ``all" is just the best performer of all the non FAQ algorithms.
Figure \ref{fig:allRelAccuracy} plots the logarithm (base 10, here and elsewhere) of the relative accuracy, i.e. $\log_{10}(\mh{f}_{FAQ}/\mh{f}_X)$, for $X \in \{$PATH, QCV, RANK, U, all$\}$.  FAQ performs significantly better than the other algorithms, outperforming all of them on $\approx 94\%$ of the problems, often by nearly an order of magnitude in terms of relative error.

\subsection{QAP Benchmark Efficiency} 
\label{sub:efficiency}

\begin{figure}[htbp]
	\centering
		\includegraphics[height=3in]{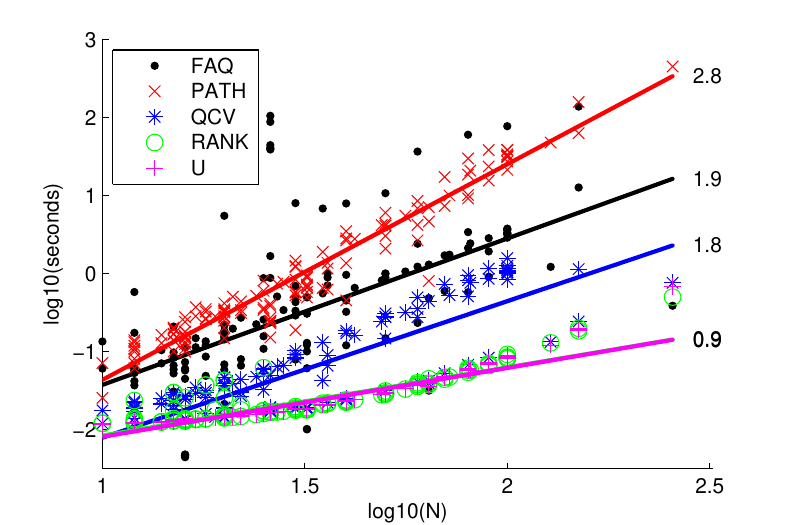}
	\caption{{\bf Absolute wall time for running each of the five algorithms on all 137 QAPLIB benchmarks}. We fit a line on this log-log plot for each algorithm; the slope is displayed beside each line. The FAQ slope is much better than the PATH slope, and worse than the others.  Note, however, the time for RANK and U appears to be superlinear on this log-log plot, suggesting that perhaps as the number of vertices increases, PATH might be faster. }
	\label{fig:allEfficiency}
\end{figure}

The utility of an approximation algorithm depends not just on its accuracy, but also its efficiency.
To empirically test these algorithms' efficiency, we compare the wall time of each of the five algorithms on all 137 QAPS in QAPLIB in Figure \ref{fig:allEfficiency}.  For each of the 5 algorithms, we fit an iteratively weighted least squares linear regression function (Matlab's ``robustfit'') to regress the logarithm of time (in seconds) onto the logarithm of the number of vertices being matched.  The numbers beside the lines indicate the slopes of the regression functions.  

The figure demonstrates that the PATH algorithm has the largest slope, significantly larger than that of FAQ.  QCV and FAQ have nearly identical slopes, which makes sense, given that the are solving very similar objective functions.  Similarly, RANK and U have very similar slopes; they are both using spectral approaches.  Note, however, that although the slope of RANK and U are smaller than that of FAQ, they both appear to be super linear in this log-log plot, suggesting that as the number of vertices increases, their compute time might exceed that of the other algorithms.  

Combined with Figure \ref{fig:allRelAccuracy}, this figure suggests that FAQ achieves the state-of-the-art trade-off between efficiency and accuracy.  
Of note is that the FAQ algorithm has a relatively high variance in wall time for these problems.  This is due to the number of Hungarian algorithms performed in the FW subroutine.  We could fix the number of Hungarian algorithms, in which case the variance would decrease dramatically.  However, in application, this variance is not problematic.

\subsection{QAP Benchmark Accuracy/Efficiency Trade-off} 
\label{sub:tradeoff}

\begin{figure}[h!]
	\centering
		\includegraphics[width=1.0\linewidth]{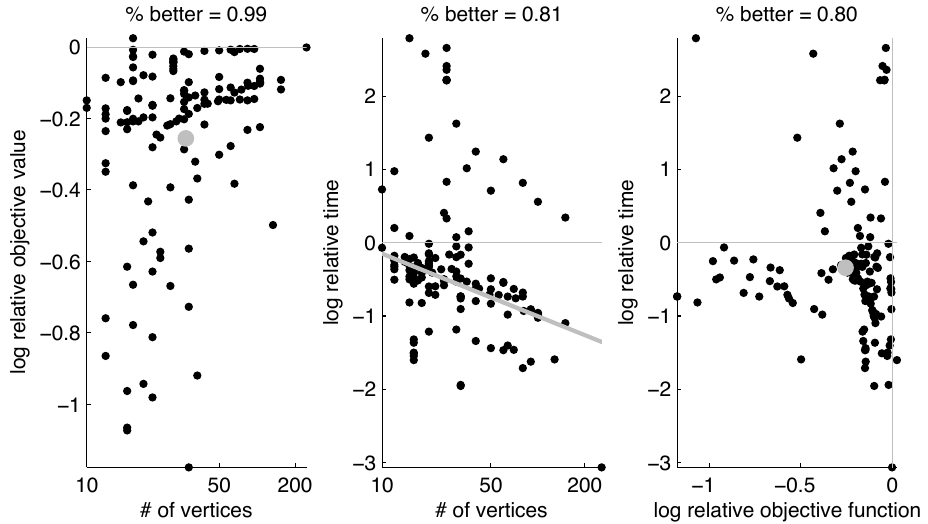}
	\caption{{\bf Comparison of FAQ with PATH in terms of both accuracy and efficiency}.  The left panel is the same as the left panel of Figure \ref{fig:allRelAccuracy}.  The middle plots the relative wall time of FAQ to PATH as a function of the number of vertices, also on a log-log scale.  The gray line is the best fit slope on this plot.  Finally, the right panel plots log relative time versus log relative objective function value, demonstrating that FAQ outperforms PATH on both dimensions on $80\%$ of the benchmarks.}
	\label{fig:tradeoff}
\end{figure}

In \cite{Zaslavskiy2009}, the authors demonstrated that PATH outperformed both QCV and U on a variety of simulated and real examples.  Figure \ref{fig:tradeoff} compares the performance of FAQ with PATH along both dimensions of performance---accuracy and efficiency---for all 137 benchmarks in the QAPLIB library.  The right panel indicates that FAQ is {\it both} more accurate and more efficient on $80\%$ of the problems (and is more accurate on $99\%$ of the benchmarks).  The middle plots the relative wall time of FAQ to PATH as a function of the number of vertices, also on a log-log scale.  The gray line is the best fit slope on this plot, suggesting that FAQ is getting exponentially faster than PATH as the number of vertices gets larger.

\subsection{QAP Directed Benchmarks}
\label{sub:directed}

Recently, Liu et al. \cite{Liu2012} proposed a modification of the PATH algorithm that adjusted PATH to be more appropriate for directed graphs
Note that our FAQ algorithm 
easily extend to directed or weighted graphs.
Liu et al. compare the performance of their algorithm (EPATH) with U, QCV, and the GRAD algorithm of \cite{Gold1996} on a set of 16 particularly difficult directed benchmarks from QAPLIB.  The EPATH algorithm achieves at least as low objective value as the other algorithms on 15 of 16 benchmarks.  Our algorithm, FAQ, performs at least as well as EPATH, U, QCV, and GRAD on all 16 benchmarks, and achieves the singular best performance on 8 of the benchmarks.  Table \ref{tab:directed} shows the numerical results comparing FAQ to EPATH and GRAD, the only algorithm considered in \cite{Liu2012} to outperform EPATH.  Note that some of the algorithms achieve the absolute minimum on 8 of the benchmark.

\begin{table}[h!]
\caption{{\bf Comparison of FAQ with optimal objective function value and previous state-of-the-art for directed graphs}.  The best (lowest) value is in \textbf{bold}. Asterisks indicate achievement of the global minimum.  The number of vertices for each problem is the number in its name (second column).}
\begin{center}
\begin{tabular}{|r|r|r||l|l|l|l|l|}
	\hline 
	          \# &  Problem &      Optimal & FAQ & EPATH & GRAD \\
	\hline 
	           1 &  lipa20a &     3683 & \textbf{3791} &     3885 &     3909 \\ 
	           2 &  lipa20b &    27076 & \textbf{27076}$^*$ &    32081 &    \textbf{27076}$^*$ \\ 
	           3 &  lipa30a &    13178 & \textbf{13571} 	&    13577 &    13668 \\ 
	           4 &  lipa30b &   151426 & \textbf{151426}$^*$ & \textbf{151426}$^*$ &   \textbf{151426}$^*$ \\ 
	           5 &  lipa40a &    31538 & \textbf{32109} 	&    32247 &    32590 \\ 
	           6 &  lipa40b &   476581 & \textbf{476581}$^*$ &   \textbf{476581}$^*$ &   \textbf{476581}$^*$ \\ 
	           7 &  lipa50a &    62093 & \textbf{62962} &    63339 &    63730 \\ 
	           8 &  lipa50b &  1210244 & \textbf{1210244}$^*$ &  \textbf{1210244}$^*$ &  \textbf{1210244}$^*$ \\ 
	           9 &  lipa60a &   107218 & \textbf{108488} &   109168 &   109809 \\ 
	          10 &  lipa60b &  2520135 & \textbf{2520135}$^*$ &  \textbf{2520135}$^*$ &  \textbf{2520135}$^*$ \\ 
	          11 &  lipa70a &   169755 & \textbf{171820} &   172200 &   173172 \\ 
	          12 &  lipa70b &  4603200 & \textbf{4603200}$^*$ &  \textbf{4603200}$^*$ &  \textbf{4603200}$^*$ \\ 
	          13 &  lipa80a &   253195 & \textbf{256073} &   256601 &   258218 \\ 
	          14 &  lipa80b &  7763962 & \textbf{7763962}$^*$ &  \textbf{7763962}$^*$ &  \textbf{7763962}$^*$ \\ 
	          15 &  lipa90a &   360630 & \textbf{363937} &   365233 &   366743 \\ 
	          16 &  lipa90b & 12490441 & \textbf{12490441}$^*$ & \textbf{12490441}$^*$ & \textbf{12490441}$^*$ \\ 
	\hline
	\end{tabular}
\end{center}
\label{tab:directed}
\end{table}

\subsection{Theoretical properties of FAQ}
\label{sub:theory}

In addition to guarantees on computational time, we have a guarantee on performance:  
\begin{prop}
If $A$ is the adjacency matrix of an asymmetric simple graph $G$, then 

$$\argmin_{D\in\mc{D}} -\text{trace}(ADA\T D\T)=\{I\}.$$
\end{prop}
{\it Proof}
Let $m$ denote the number of edges in $G$.  As $G$ is asymmetric,

$$\text{trace}(APA\T P\T)<2m=\text{trace}(AA\T )$$ for any $P\neq I.$  By the Birkhoff-von Neuman Theorem, ${\mathcal D}$ is the convex hull
of $\mc{P}$, i.e., for all $D \in {\mathcal D}$, there exists
constants $\{ \alpha_{D,P} \}_{P \in \mc{P}}$ such that
$D=\sum_{P \in \mc{P}}\alpha_{D,P}P$ and $\sum_{P \in \mc{P}}\alpha_{D,P}=1$.
Thus, if $D$ is not the identity matrix,

\begin{align*}
\text{trace}(ADA\T D\T)&=\sum_{P}\sum_{Q}\alpha_{D,P}\alpha_{D,Q}\text{trace}(APA\T Q\T)\\
&=\sum_{P}\alpha_{D,P}^2\underbrace{\text{trace}(APA\T P\T)}_{<2m\text{ if }P\neq I}+\sum_{P}\sum_{Q\neq P}\alpha_{D,P}\alpha_{D,Q}\underbrace{\text{trace}(APA\T Q\T)}_{\leq 2m}\\
&<2m
\end{align*}
as $D\neq I$ implies that $\alpha_{D,P}>0$ for some $P\neq I$. $\blacksquare$

\begin{rem}
	\emph{Note that it trivially follows from Proposition 1 that if $A$ and $B$ are the adjacency matrices of asymmetric isomorphic simple graphs, then the minimum objective function value of rGMP is equal to the minimum objective function value of GMP.}
\end{rem}

\begin{rem}
\emph{ For the convex quadratic GM relaxation, Eq. (\ref{eq:GMc}),
in general } 
$\argmin_{D\in\mc{D}} \|AD-DA\|^2_F\neq\{I\},$
\emph{even if $G$ (the graph with adjacency matrix $A$) is asymmetric.  Indeed, degree regular graphs provide a simple counterexample, as in this case} 
$J\in\argmin_{D\in\mc{D}} \|AD-DA\|^2_F.$
\emph{We will empirically explore the ramifications of this phenomena further in Section \ref{sub:connectome}}
\end{rem}
% {\it Proof}
% Let $A=PBP\T$, so that $\langle A, PBP\T\rangle=2m$, where $m$ is the number of edges in $A$.  By way of contradiction, assume that there exists a $D\in\mc{D}\backslash\mc{P}$ such that $\langle A, DBD\T\rangle < \langle A, PBP\T\rangle = 2m$.
% %(remember that we are minimizing the negative Euclidean inner product)
% It follows that $(DBD\T)_{uv}> 1$ for some $(u,v)$.  However, the submultiplicativity of the norm induced by the $\ell_{\infty}$ norm yields
% $\norm{Dx}_\infty \leq \norm{D}_{\infty,\infty} \norm{x}_\infty$.  Applying this twice (once for each doubly stochastic matrix multiplication) yields the desired contradiction.
% % Consider $d_i=\langle D, \text{col}_i(BD\T) \rangle$, where $\text{col}_i(\cdot)$ indicates the $i^{th}$ column of the matrix.  $d_i \leq 1$ for all $i \in [n]$, therefore, our result holds.
% $\blacksquare$

This result says that, when solving the GI problem, nothing is lost by relaxing the indefinite GM problem as done by FAQ.  Recent results also show this is almost surely the case (in a broad class of random graphs) when relaxing the indefinite GM problem, even in the non-isomorphic graph setting \cite{lyzinski2014graph} (and is again almost surely {\it not} the case when relaxing the convex GM formulation).  Combined, this serve to provide theoretical justification for our FAQ procedure. 

\subsection{Multiple Restarts} 
\label{sub:multiple_restarts}

Due to the indefiniteness of the relaxation of (\ref{eq:GM}), the feasible region may be rife with local minima.  
As a result, our FAQ procedure is sensitive to the chosen initial position.    
With this in mind, we propose a variant of the FAQ algorithm in which we run the FAQ procedure with
multiple initializations.  The algorithm outputs the best FAQ iterate over all the initializations.
For each initialization, we  sample $K \in \mc{D}$, a random doubly stochastic matrix, using 10 iterations of Sinkhorn balancing \cite{Sinkhorn1964}, and let our initialization be $P^{(0)}=(J+K)/2$, where $J$ is the doubly stochastic barycenter.  Fixing the number of restarts, this variant of FAQ still has $O(n^3)$ complexity.    

Table \ref{tab:restarts} shows the performance of running FAQ 3 and 100 times, reporting only the best result (indicated by FAQ$_3$ and FAQ$_{100}$, respectively), and comparing it to the best performing result of the five algorithms 
(including the original FAQ)---note that the best performing of the original five tested was always FAQ. Note that we only consider the 16 particularly difficult benchmarks for this evaluation. FAQ only required three restarts to outperform all other approximate algorithms on all 16 of 16 difficult benchmarks.  Moreover, after 100 restarts, FAQ finds the absolute minimum on 3 of the 16 benchmarks.  Note that no other algorithm ever achieved the absolute minimum on any of these benchmarks. 

\begin{rem}
\emph{Another natural starting position for FAQ is the doubly stochastic solution of rGMP, the relaxed Equation (\ref{eq:GMc}).  Promising results in this direction are pursued further in \cite{lyzinski2014graph}.}
\end{rem}

\begin{table}[h!]
\caption{{\bf Comparison of FAQ with optimal objective function value and the best result  on the undirected benchmarks}.  Note that FAQ restarted 100 times finds the optimal objective function value in 3 of 16 benchmarks, and that FAQ restarted 3 times finds a minimum better than the previous state-of-the-art on all 16 particularly difficult benchmarks.}
\begin{center}
\begin{tabular}{|r|r|r||l|l|l|l|l|}
\hline
\# & Problem  &   Optimal    & FAQ$_{100}$ & FAQ$_{3}$ & FAQ$_1$\\
\hline
1&    chr12c &   11156 &    \textbf{12176} &   13072 & 13072 \\
2&    chr15a &    9896 &    \textbf{9896}$^*$ &   17272 &  19086 \\
3&    chr15c &    9504 &    \textbf{10960} &   14274 &  16206 \\
4&   chr20b &    2298 &     \textbf{2786} &    3068 &    3068 \\
5&    chr22b &    6194 &    \textbf{7218} &    7876 &   8482 \\
6&    esc16b & 	292 & 		\textbf{292}$^*$ & 294 &    296 \\
7& 	   rou12 &  235528 &  \textbf{235528}$^*$ &  238134 &    253684 \\
8& 	   rou15 &  354210 &  \textbf{356654} &  371458 &    371458 \\
9&      rou20 &  725522 &  \textbf{730614} &  743884 &    743884 \\
10&    tai10a &  135028 &  \textbf{135828} &  148970 &    152534 \\
11&    tai15a &  388214 &  \textbf{391522} &  397376 &    397376 \\
12&    tai17a &  491812 &  \textbf{496598} &  511574 &    529134 \\
13&    tai20a &  703482 &  \textbf{711840} &  721540 &    734276 \\
14&    tai30a & 1818146 & \textbf{1844636} & 1890738 &  1894640 \\
15&    tai35a & 2422002 & \textbf{2454292} & 2460940 &  2460940 \\
16&    tai40a & 3139370 & \textbf{3187738} & 3194826 &  3227612 \\
    \hline
\end{tabular}
\end{center}
\label{tab:restarts}
\end{table}

\subsection{Brain-Graph Matching} 
\label{sub:connectome}

The \emph{Caenorhabditis elegans} (\emph{C. elegans}) is a small worm (nematode) with $302$ labeled neurons (in the hermaphroditic sex).  The chemical connectome of \emph{C. elegans} is a weighted directed graph on 279 vertices, with edge weights in $\{0,1,2,\ldots\}$ counting the number of directed chemical synapses between the neurons \cite{WhiteBrenner86, Varshney2011}.    
We conduct the following synthetic experiments.  
For $A=(A_{uv})$ count the number of synapses (in $\{0,1,2,\ldots\}$) from neuron $u$ to neuron $v$ in the C. elegans connectome, for all $u,v$.  For $k =1,2,\ldots, 1000$, we choose $Q_{(k)}$ uniformly at random from $\mc{P}$, and let $B_{(k)}=Q_{(k)} A {Q_{(k)}}\T$, effectively shuffling the vertex labels of the connectome.  Then, we match graphs $A$ to $B_{(k)}$.  We define accuracy as the fraction of vertices correctly assigned (i.e. unshuffled).

\begin{figure}[h!]
	\centering
		  \includegraphics[width=0.7\linewidth]{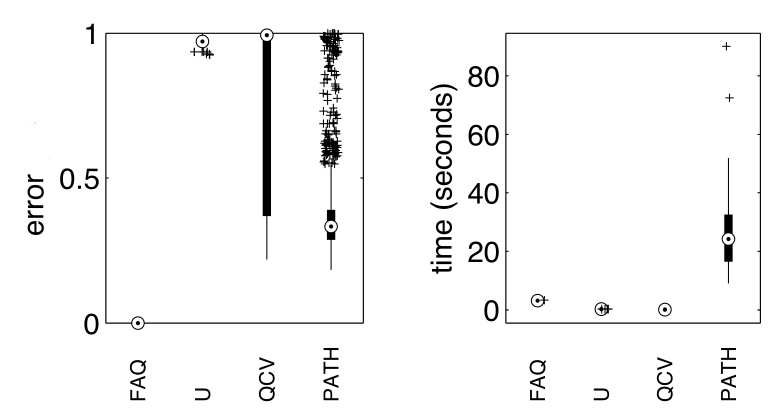}
	\caption{{\bf Performance of U, QCV, PATH, and FAQ on synthetic C. elegans connectome data, graph matching the true chemical connectome with permuted versions of itself}.  Error is the fraction of vertices incorrectly matched.  Circle indicates the median, thick black bars indicate the quartiles, thin black lines indicate extreme but non-outlier points, and plus signs are outliers. The left panel indicate error (fraction of misassigned vertices), and the right panel indicates wall time on a 2.2 GHz Apple MacBook.   FAQ always obtained the optimal solution, whereas none of the other algorithms ever found the optimal.  FAQ also ran very quickly, nearly as quickly as U and QCV, and much faster than PATH, even though the FAQ implementation is in Matlab, and the others are in C.}
	\label{fig:connectomes}
\end{figure}

Figure \ref{fig:connectomes} displays the results of FAQ (initialized at $J$) along with U, QCV,  and PATH.  The left panel indicates that FAQ \emph{always} perfectly unshuffles the chemical connectome, whereas none of the other algorithms ever perfectly unshuffles the graph.  
In light of Proposition 1, this is unsurprising.  Indeed, there is a unique automorphism for this connectome, and the graph is asymmetric.  For any choice of $Q_{(k)}$, the indefinite problem (\ref{eq:GM}) therefore has a unique solution, namely $Q_{(k)}\T$.  FAQ finds this global minima in all the cases.  Contrast this with equation (\ref{eq:GMc})---the objective function of PATH and QCV---which could have multiple global minima in $\mc{D}$.  This could account for the high variance in the performance of QCV and PATH in Figure \ref{fig:connectomes}. 
 
The right panel compares the wall time of the various algorithms, running on an 2.2 GHz Apple MacBook. Note that we only have a Matlab implementation of FAQ, whereas the other algorithms are implemented in C.  
Unlike in the QAPLIB benchmarks, FAQ runs nearly as quickly as both U and QCV; and as expected, FAQ runs significantly faster than PATH.  
We also ran FAQ on a binarized symmeterized versions of the chemical connectome graph $\tilde A$ (i.e. $\tilde A_{uv}=1$ if and only if $A_{uv}\geq 1$ or $A_{vu} \geq 1$).  The resulting errors are nearly identical to those presented in Figure \ref{fig:connectomes}, although speed increased for FAQ by more than a factor of two. Note that the number of vertices in this brain-graph matching problem---279---is larger than the largest of the 137 benchmarks used above.

\section{Discussion}
\label{sec:discussion}

This work presents the FAQ algorithm, a fast algorithm for approximately matching very large graphs.  Our key insight was to continuously relax the {\it indefinite} formulation of the GMP in the FW implementation.
We demonstrated that not only is FAQ cubic in time, but also its leading constant is quite small---$10^{-9}$---suggesting that it can be used for graphs with hundreds or thousands of vertices (\S \ref{sub:const}).  

Moreover, FAQ achieves better accuracy than previous state-of-the-art approximate algorithms on on over $93\%$ of the 137 QAPLIB benchmarks (\S \ref{sub:qap_benchmarks}), is faster than the state-of-the-art PATH algorithm (\S \ref{sub:efficiency}), and is both faster and achieves at least as low performance as PATH on over $80\%$ of the tested benchmarks (\S \ref{sub:tradeoff}),  including both directed and undirected graph matching problems (\S \ref{sub:directed}).  
In addition to the theoretical guarantees of cubic run time, we provide theoretical justification for optimizing the indefinite GM formulation, Equation (\ref{eq:GM}) as opposed to Equation (\ref{eq:GMc}) (\S \ref{sub:theory}).  Indeed, the indefinite formulation (and {\it not} the convex formulation) has the property that when matching asymmetric isomorphic graphs, the unique global minimum of the indefinite relaxation is the isomorphism between the two graphs.  

Because rQAP is indefinite, we also considered FAQ with multiple restarts, and achieve improved performance for the particularly difficult benchmarks using only three restarts (\S \ref{sub:multiple_restarts}).    
Finally, we used FAQ to match permuted versions of the C. elegans connectomes (\S \ref{sub:connectome}). Of the four state-of-the-art algorithms considered, FAQ achieved perfect performance $100\%$ of the time, whereas none of the other three algorithms ever matched the connectomes perfectly.  Moreover, FAQ ran comparably fast to U and QCV and significantly faster than PATH, even though FAQ is implemented in Matlab, and the others are implemented in C.  Note that these connectomes have 279 vertices, more vertices than the largest QAP benchmarks.

\subsection{Related Work}

Our approach is quite similar to other approaches that have recently appeared in the literature.  Perhaps its closest cousins include \cite{Zaslavskiy2009, Zaslavskiy2010} and \cite{Escolano2011}, which are all of the ``PATH'' following variety.  These algorithms begin by relaxing the convex objective function in (\ref{eq:GMc}), while FAQ begins by relaxing the indefinite objective function in (\ref{eq:GM}). 
Although the convex relaxation is efficiently solvable, the obtained solution is almost surely incorrect (for a broad class of random graphs)
 and the correct solution is often not obtained even post projection \cite{lyzinski2014graph}.  The indefinite relaxation however, almost surely yields the correct solution when exactly solved (for a broad class of random graphs) 
 \cite{lyzinski2014graph}.  %Zaslavskiy et al. considers FAQ but discards   \cite{Zaslavskiy2009} because they did not like projecting onto the set of permutations matrices. 
With this in mind, it is unsurprising that FAQ outperforms PATH on nearly all benchmark problems.  

Others have considered similar relaxations to PATH, but usually in the context of finding lower bounds  \cite{Anstreicher2001} or as subroutines for finding exact solutions \cite{Brixius2000}.  Our work seems to be the first to utilize the precise algorithm described in Algorithm \ref{alg:1} to find fast approximate solutions to GMP and QAP.

\subsection{Future Work}

Even with the very small lead order constant for FAQ, as $n$ increases, the computational burden gets quite high.  For example, extrapolating the curve of Figure \ref{fig:scaling}, this algorithm would take about 20 years to finish (on a standard laptop) when $n=100,000$.  We hope to be able to approximately solve rQAP on graphs much larger than that, given that the number of neurons even in a fly brain, for example, is $\approx 250,000$.  More efficient algorithms and/or implementations are required for such massive graph matching. Although a few other state-of-the-art algorithms were more efficient than FAQ, their accuracy was significantly worse.  We are actively working on combining FAQ with dimensionality reduction procedures to achieve the desired scaling from FAQ \cite{Lyzinski2013}.

% So the search continues to find approximate graph matching algorithms with scaling rules like QCV, U or RANK, but performance like FAQ.

We are also pursuing additional future work to generalize FAQ in a number of ways:  
\begin{itemize}
\item In addition to the theoretical results contained in Section \ref{sub:theory}, we have studied the properties of the convex and indefinite GMP relaxations in a very general random graph model \cite{lyzinski2014graph}.  
Under some general assumptions on the random graph model, the indefinite relaxation of (\ref{eq:GM}), and {\it not} the convex relaxation of (\ref{eq:GMc}) is the provably correct approach \cite{lyzinski2014graph}.
\item Many (brain-) graphs of interest will be errorfully observed \cite{Priebe2011}, that is, vertices might be missing and putative edges might exhibit both false positives and negatives.  Explicitly dealing with this error source is both theoretically and practically of interest \cite{VP11_unlabeled}.  
\item For many brain-graph matching problems, the number of vertices will not be the same across the brains.  Recent work from \cite{Zaslavskiy2009, Zaslavskiy2010} and \cite{Escolano2011} suggest that extensions in this direction would be both relatively straightforward and effective. 
\item Often, a partial matching of the vertices is known a priori, and we can modify FAQ to leverage these seeded vertices to match the remaining unseeded vertices \cite{sgm2}.
\item The most costly subroutine in FAQ is solving LAPs.  Fortunately, the LAP is a linear optimization problem with linear constraints.  As a result, a number of parallelized optimization strategies could be implemented on this problem \cite{Boyd04a}. 
\item Often, real data adjacency matrices have certain special properties, namely sparsity, which makes faster LAP subroutines \cite{Jonker1987} and more efficient algorithms (such as ``active set'' algorithms) readily available for further speed increases.
\item In many graph settings, we have some prior information that could easily be incorporated into the GM problem in the form of vertex attributes and features.  For example, in brain graphs we know the position of the vertex in the brain, the vertex's cell type, etc.  These could be used to measure a ``dissimilarity'' between vertices and are easily incorporated into FAQ's objective function to better match the graphs.  
\item We are working to extend FAQ to match multiple graphs simultaneously.
\end{itemize}

\subsection{Concluding Thoughts}

In conclusion, this manuscript has presented the FAQ algorithm for approximately solving the quadratic assignment problem.  FAQ is theoretically justified, fast, effective, and easily generalizable.  Our algorithm achieves state-of-the-art matching performance and efficiency on a host of benchmark QAP problems and connectome data sets.  Yet, the $\mc{O}(n^3)$ complexity remains too slow to solve many problems of interest, and issues of scalability need be addressed.  To facilitate further development and applications, all the code and data used in this manuscript is available from the first author's website, \url{http://jovo.me}.  We have further incorporated FAQ (as sgm.R) into the open-source R package, igraph, available for download at \url{https://github.com/igraph/xdata-igraph/}.  MATLAB code is also available at \url{https://github.com/jovo/FastApproximateQAP/tree/master/code/FAQ}.

\section*{Acknowledgment}

The authors would like to acknowledge 
Lav Varshney for providing the data.

\end{document}